\documentclass[11pt]{amsart}
\usepackage{amsmath,amssymb,amsfonts,url,mathptmx}
\numberwithin{equation}{section}

\newtheorem{thm}{Theorem}[section]
\newtheorem{lem}{Lemma}[section]
\newtheorem{rem}{Remark}[section]
\newtheorem{prop}{Proposition}[section]
\newtheorem{cor}{Corollary}[section]
\newcommand{\hdot}{^\text{\r{}}\hspace{-.33cm}H}

\begin{document}
\title[Singular Liouville System]{Degree counting theorems for singular Liouville systems} 

\author{Yi Gu}
\address{Department of Mathematics\\
        University of Florida\\
        1400 Stadium Rd\\
        Gainesville FL 32611}
\email{yigu57@ufl.edu}
\author{Lei Zhang}\footnote{Lei Zhang is partially supported by a Simons Foundation Collaboration Grant}
\address{Department of Mathematics\\
        University of Florida\\
        1400 Stadium Rd\\
        Gainesville FL 32611}
\email{leizhang@ufl.edu}

\date{\today}

\begin{abstract} Let $(M,g)$ be a compact Riemann surface with no boundary and $u=(u_1,...,u_n)$ be a solution of the following singular Liouville system:
\begin{equation*}
\Delta_g u_i+\sum_{j=1}^na_{ij}\rho_j(\frac{h_je^{u_j}}{\int_M h_j e^{u_j}dV_g}-\frac{1}{vol_g(M)})=\sum_{t=1}^N4\pi \gamma_t( \delta_{p_t}-\frac{1}{vol_g(M)}),
\end{equation*}
where $i=1,...,n$,
 $h_1,...,h_n$ are positive smooth functions, $p_1,...,p_N$ are distinct points on $M$, $\delta_{p_t}$ are Dirac masses, $\rho=(\rho_1,...,\rho_n)$ ($\rho_i\ge 0)$ and $(\gamma_{1},...,\gamma_{N})$ ($\gamma_{t}>-1$ ) are constant vectors. If the coefficient matrix $A=(a_{ij})_{n\times n}$ satisfies standard assumptions we identify a family of critical hyper-surfaces $\Gamma_k$ for $\rho=(\rho_1,..,\rho_n)$ so that a priori estimate of $u$ holds if $\rho$ is not on any of the $\Gamma_k$s. Thanks to the a priori estimate, a topological degree for $u$ is well defined for $\rho $ staying between every two consecutive $\Gamma_k$s. In this article we establish this degree counting formula which depends only on the Euler Characteristic of $M$ and the location of $\rho$. Finally if the Liouville system is defined on a bounded domain in $\mathbb R^2$ with Dirichlet boundary condition, a similar degree counting formula that depends only on the topology of the domain and the location of $\rho$ is also determined.
\end{abstract}


\maketitle

\section{Introduction}
In this article we study the following Liouville system defined on a compact Riemann surface $(M,g)$ with no boundary:
\begin{eqnarray}
\Delta_gu_i^*+\sum_{j=1}^n\rho_j a_{ij}(\frac{h_j^*e^{u_j^*}}{\int_M h_j^* e^{u_j^*}dV_g}-\frac{1}{vol_g(M)})=\sum_{l=1}^N 4\pi \gamma_{il}(\delta_{p_l}-\frac{1}{vol_g(M)}), \nonumber \\
\mbox{ for } \quad i\in  I:=\{1,...,n\}, \quad \gamma_{il}>-1, \mbox{ for } \, i\in I, l=1,...,N;\label{liou-sys}
\end{eqnarray}
where $h_1^*,...,h_n^*$ are positive smooth functions on $M$, $\rho_1,...,\rho_n$ are nonnegative constants, $vol_g(M)$ is the volume of $M$, $p_1,...,p_N$ are distinct points on $M$, $\delta_{p_l}$ are singular sources at $p_l$ and $\gamma_{il}>-1$ ($i=1,..,n$, $l=1,...N$) are constants as well. Equation (\ref{liou-sys}) is called Liouville system if all the entrees in the coefficient matrix $A=(a_{ij})_{n\times n}$ are nonnegative.

System (\ref{liou-sys}), in its generality, covers a large number of models in different subjects of mathematics, physics and other disciplines as well. In physics Liouville systems can be derived from the mean field limit of point vortices of the Euler flow ( see \cite{caglioti-1,caglioti-2,kiess,cha-kiess}). The study of Liouville systems finds applications in nonabelian Chern-Simons-Higgs theory (\cite{dunne, hong,jackiw,yang}) and the electroweak theory (see \cite{ambjorn,nolasco-ta-1,nolasco-ta-2,nolasco-ta-3,spruck-yang-1,spruck-yang-2,taran,yang}).Various Liouville systems are also used to describe models in theories of chemotaxis (\cite{childress,keller}), the physics of charged particle beams \cite{bennet,debye,kiess-2,kiess-le}, and other gauge field models \cite{dz,kim}. Even if the system is reduced to one equation, it has profound background in geometry: if the equation has no singular source, it interprets the Nirenberg problem of prescribing Gauss curvature; if the equation has singular sources, the solution represents a metric with conic singularity \cite{kuo-lin}.  It is just impossible to overestimate the importance of Liouville systems.

One of the main goals in the study of Liouville system is to identity the role that the topology of $M$ plays in the structure of solutions. In particular, people seek to identity a family of hyper-surfaces for $\rho:=(\rho_1,...,\rho_n)$, so that if $\rho$ does not belong to these hyper-surfaces, a priori estimate of $u$ holds and the Leray-Schauder degree can be defined. The explicit computation of the Leray-Schauder degree, which depends on the topology of $M$, gives rise to existence of solution if the degree is not zero.  Usually the identification of critical hyper-surfaces requires detailed study of blowup solutions, and it is well known that local, geometric information, such as the Gauss curvature plays a crucial role in determining the asymptotic behavior of blowup solutions, the main purpose of this article is to establish a link between local analysis, the structure of solutions and the topology of 2-manifolds for a class of singular Liouville systems.

If the system is reduced to Liouville equation, Chen and Lin completed the program in a series of pioneering works \cite{chenlin1,chenlin2,chenlin3}. The readers may read into \cite{licmp,li-shafrir, zhangcmp, bclt, bart2, zhangccm} for background and related discussions. Chen-Lin's work was extended by Lin and the second author \cite{lin-zhang-1,linzhang2,lin-zhang-3} to Liouville systems with no Dirac sources. Since singular sources have significant geometric applications, the main purpose of this article is to extend Lin-Zhang's degree counting formula to systems with Dirac poles.

For the coefficient matrix $A$ we postulate two conditions: The first one is called a standard assumption:
\begin{equation}\nonumber
(H1):\quad A \mbox{ is symmetric, non-negative, irreducible and invertible}.
\end{equation}
Here we note that $A$ being irreducible means there is no partition of the index set $I:=\{1,...,n\}$ into two disjoint subsets $I=I_1\cup I_2$ such that $a_{ij}=0$ for all $i\in I_1$ and $j\in I_2$. In other words the Liouville system cannot be written as two separated sub-systems.
The second assumption, which is made on the inverse of $A^{-1}=(a^{ij})_{n\times n}$, is called a strong interaction assumption: For $I=\{1,...,n\}$,
$$
(H2):\left\{ \begin{array}{ll}
a^{ii}\le 0,\,\, \,\, \forall i\in I,\qquad a^{ij}\ge 0,\quad \forall i\neq j, i,j\in I,\\  \\
\sum_{j\in I}a^{ij}\ge 0, \,\, \forall i\in I.
\end{array}
\right.
$$
The reason that  $(H2)$ is called a strong interaction assumption can be justified from the following two examples: For $n=2$, the matrix
$$A=\left(\begin{array}{cc}
a_{11} & a_{12} \\
a_{12} & a_{22}
\end{array}
\right)
$$
satisfies $(H1)$ and $(H2)$ if and only if $a_{ij}\ge 0$, $\max(a_{11},a_{22})\le a_{12}$, and $det(A)\neq 0$. For $n=3$, the following matrix
$$A_1=\left(\begin{array}{ccc}
0 & a_1 & a_2 \\
a_1 & 0 & a_3\\
a_2 & a_3 & 0
\end{array}
\right )
$$
satisfies both $(H1)$ and $(H2)$ if and only if $a_i>0$ and $a_i+a_j\ge a_k$ for $i,j,k$ all different from one another.

\medskip

The second main assumption is that around each singular source, the strength of the singular source for each component is the same: $\gamma_{il}=\gamma_l>-1$ for all $i=1,...,n$. This assumption is crucial to for ruling out all partial blowups later. 
Also for convenience we assume that the volume of the manifold is $1$, thus (\ref{liou-sys}) can be written as
\begin{equation}\label{syst-2}
\Delta_gu_i^*+\sum_{j=1}^n\rho_j a_{ij}(\frac{h_j^*e^{u_j^*}}{\int_M h_j^* e^{u_j^*}dV_g}-1)=\sum_{l=1}^N 4\pi \gamma_{l}(\delta_{p_l}-1)
\end{equation}

Around each singular source, the leading term of $u_i^*$ is a logarithmic function that comes from the following Green's function $G(x,q)$:

\begin{equation}\label{green-f}
\left\{\begin{array}{ll}
-\Delta_xG(x,q)=\delta_q-1,\\
\\
\int_MG(x, q)dx=0.
\end{array}
\right.
\end{equation}
It is a common practice to define
$$u_i=u_i^*-4\pi\sum_{l=1}^N\gamma_lG(x,p_l),$$
and rewrite (\ref{syst-2}) as
\begin{equation}\label{main-sys-2}
\Delta_gu_i+\sum_{j=1}^na_{ij}\rho_j(\frac{h_je^{u_j}}{\int_Mh_je^{u_j}}-1)=0,\quad i=1,...,n,
\end{equation}
where
$$h_i(x)=h_i^*(x)exp\{-\sum_{l=1}^N4\pi \gamma_lG(x,p_l)\}, $$
which implies that around
each singular source, say, $p_l$, in local coordinates, $h_j$ can be written as
$$h_j(x)=|x|^{2\gamma_l}g_j(x)$$
for some positive, smooth function $g_j(x)$.

Obviously, equation (\ref{main-sys-2}) remains the same if $u_i$ is replaced by $u_i+c_i$ for any constant $c_i$. Thus we might assume that each component of $u=(u_1,...,u_n)$ is in
$$ \hdot^{1}(M):=\{v\in L^2(M);\quad \nabla v\in L^2(M), \mbox{and }\,\, \int_M v dV_g=0\}. $$
Then equation (\ref{main-sys-2}) is the Euler-Lagrange equation for the following nonlinear functional $J_{\rho}(u)$ in $\mathring{H}^{1}(M)$:
$$J_{\rho}(u)=\frac 12\int_M\sum_{i,j=1}^na^{ij}\nabla_gu_i\nabla_gu_jdV_g-{\sum_{i=1}^n\rho_i\log \int_M h_ie^{u_i}dV_g}.$$

Let $\mathbb N^+$ be the set of positive integers. We shall use the following notation:
$$\Sigma:=\{8m\pi+\sum_{p_l\in A} 8\pi(1+\gamma_l);\quad A\subset \{p_1,...,p_N\},\quad m\in \mathbb N^+\cup\{0\} \quad \}\setminus \{0\}. $$
Writing $\Sigma$ as
\begin{equation}\label{e-nk}
\Sigma=\{8\pi n_k \quad | \quad n_1<n_2<... \quad \}
\end{equation}
we first establish the following a priori estimate:
\begin{thm}\label{a-pri-main}
Let $A=(a_{ij})_{n\times n}$ satisfy $(H1)$ and $(H2)$. For $k\in \mathbb N^+\cup \{0\}$, and
\begin{eqnarray*}
\mathcal{O}_k=\{(\rho_1,...,\rho_n)| \,\,  \rho_i\ge 0, i\in I;\quad \mbox{ and } \\
8\pi n_k\sum_{i\in I} \rho_i<\sum_{i,j\in I}a_{ij}\rho_i\rho_j<8\pi n_{k+1}\sum_{i\in I} \rho_i. \quad \}
\end{eqnarray*}
Suppose $h_i^*$s are positive and $C^1$ functions on $M$ and $K$ is a compact subset of $\mathcal{O}_k$. Then there exists a constant $C$ such that for any solution
$u=(u_1,...,u_n)$ of (\ref{main-sys-2}) with $\rho\in K$ and $u_i\in \,\hdot^1(M)$, we have
$$|u_i(x)|\le C, \quad \mbox{ for } i\in I, \quad \mbox{ and } \quad x\in M. $$
\end{thm}

Note that the set $\mathcal{O}_k$ is bounded if all $a_{ii}>0$ and is unbounded if $a_{ii}=0$ for some $i$. By Theorem \ref{a-pri-main}, the critical parameter set for (\ref{main-sys-2}) is
$$\Gamma_k=\{\rho;\quad 8\pi n_k\sum_{i\in I}\rho_i=\sum_{i,j\in I} a_{ij}\rho_i\rho_j \}. $$
Thanks to Theorem \ref{a-pri-main}, for $\rho\not \in \Gamma_k$, we can define the nonlinear map $T_{\rho}=(T^1,...,T^n)$ from \,\, $\hdot^{1,n}=\hdot^1(M)\times...\times \hdot^1(M)$ to \,\, $\hdot^{1,n}$ by
$$T^i=-\Delta_g^{-1}(\sum_{j\in I}a_{ij}\rho_j(\frac{h_je^{u_j}}{\int_M h_je^{u_j}}-1)),\quad i\in I. $$

Obviously $T_{\rho}$ is compact from \,$\hdot^{1,n}$ to itself. Then we can define the Leray-Schauder degree of (\ref{main-sys-2}) by
$$d_g=deg(I-T_{\rho}; B_R, 0), $$
where $R$ is sufficiently large and $B_R=\{u;\quad u\in \hdot^{1,n},\,\, \mbox{ and } \sum_{i=1}^n \|u_i\|_{H^1}< R\}$.  By the homotopic invariance and Theorem \ref{a-pri-main}, $d_{\rho}$ is constant for $\rho\in \mathcal{O}_k$ and is independent of $h=(h_1,...,h_n)$.

To state our degree counting formula for $d_{\rho}$ we consider the following generating function $g$:
$$g(x)=(1+x+x^2+...)^{-\chi(M)+N}\Pi_{l=1}^N(1-x^{1+\gamma_l}), $$
where $\chi(M)=2-2g_e(M)$ is the Euler Characteristic of $M$ ($g_e(M)$ is the genus of $M$). It is obvious to observe that if $-\chi(M)+N>0$,
$$(1+x+x^2+...)^{-\chi(M)+N}=(1-x)^{\chi(M)-N}. $$
Writing $g(x)$ in the following form
$$g(x)=1+b_1x^{n_1}+b_2x^{n_2}+...  ,$$
we use $b_1,b_2,...$ to describe our degree counting theorem:

\begin{thm}\label{main-thm}
Let $d_{\rho}$ be the Leray-Schauder degree for (\ref{main-sys-2}). Suppose
$$8\pi n_k\sum_{i=1}^n\rho_i<\sum_{i,j}a_{ij}\rho_i\rho_j<8\pi n_{k+1}\sum_{i=1}^n\rho_i, $$
then
$$d_{\rho}=\sum_{j=0}^k b_j,\quad \mbox{ where } \quad b_0=1. $$
\end{thm}
For most applications $\gamma_l$ are positive integers, which implies that
$$\Sigma=\{8\pi m ;\quad m\in \mathbb N^+\}.$$
 Thus in this case ($\gamma_l\in \mathbb N^+$) if $\chi(M)\le 0$ we have
\begin{eqnarray}\label{formula-g}
g(x)=(1+x+x^2+...)^{-\chi(M)}\Pi_{l=1}^N\frac{1-x^{1+\gamma_l}}{1-x}\\
=(1+x+x^2+...)^{-\chi(M)}\Pi_{l=1}^N(1+x+...+x^{\gamma_l})\nonumber \\
=1+b_1x+b_2x^2+...+b_kx^k+.. \nonumber
\end{eqnarray}
obviously $b_j\ge 0$ for all $j\ge 1$, which implies
$$d_{\rho}=1+\sum_{j=1}^k b_j>0. $$

\begin{cor}\label{cor-1}
Suppose all $\gamma_l\in \mathbb N^+$ and $\chi(M)\le 0$. Then $d_{\rho}>0$ if
$$\sum_{ij\in I}a_{ij}\rho_i\rho_j\neq 8\pi m\sum_{i\in I} \rho_i \quad \forall m\in \mathbb N^+. $$
Thus (\ref{main-sys-2}) always has a solution in this case.
\end{cor}

For an open, bounded smooth domain in $\mathbb R^2$, we are also interested in the following system of equations:
\begin{equation}\label{liou-loc}
\left\{\begin{array}{ll}\Delta u_i+\sum_{j=1}^na_{ij}\rho_j\frac{h_j^*e^{u_j}}{\int_{\Omega}h_j^*e^{u_j}}=4\pi \sum_{l=1}^N\gamma_l\delta_{p_l},\quad i\in I, \\
\\
u_i|_{\partial \Omega}=0, \quad i\in I,
\end{array}
\right.
\end{equation}
where $h_1^*$,...,$h_n^*$ are smooth functions on $\bar \Omega$ and $p_1,...,p_N$ are distinct points in the interior of $\Omega$.

Let
$$g(x)=(1+x+x^2+...)^{-\chi(\Omega)+N}\Pi_{l=1}^N(1-x^{1+\gamma_l})=\sum_{j=0}^{\infty}b_jx^{n_j}, $$
where $\chi(\Omega)=1-g_e(\Omega)$ ($g_e(\Omega)$ is the number of holes bounded by $\Omega$) is the Euler Characteristic number of $\Omega$, and $b_0=1$. Then we have

\begin{thm}\label{thm-local}
Suppose
$$8\pi n_k \sum_{i\in I}\rho_i<\sum_{i,j\in I} a_{ij}\rho_i\rho_j<8\pi n_{k+1} \sum_{i\in I} \rho_i. $$
Then $d_{\rho}=\sum_{j=0}^kb_j$. If $\gamma_1,...,\gamma_N\in \mathbb N^+$, $\Omega$ is not simply connected and $\sum_{ij}a_{ij}\rho_i\rho_j\neq 8\pi m\sum_i \rho_i$
for any $m\in \mathbb N$, we have $d_{\rho}>0$ and the existence of a solution to (\ref{liou-loc}).
\end{thm}

If the Liouville system on $(M,g)$ is written as
\begin{equation}\label{main-sys-5}
\Delta_g u_i^*+\sum_{j\in I} a_{ij} h_j^*e^{u_j^*}=4\pi \sum_{l=1}^N \gamma_l \delta_{p_l}, \quad i\in I,
\end{equation}
with the same assumptions on $A$, $h_i^*$, $\gamma_l$ and $vol(M)=1$,
we first remark that (\ref{main-sys-5}) is a special case of (\ref{main-sys-2}). Indeed,
integrating (\ref{main-sys-5}) on both sides, we have
$$\sum_{j\in I} a_{ij} \int_M h_j^*e^{u_j^*}=4\pi \sum_l \gamma_l. $$
Thus
\begin{equation}\label{ener-8}
\int_Mh_i^*e^{u_i^*}=4\pi \sum_{j\in I} a^{ij}(\sum_l \gamma_l),\quad i\in I.
\end{equation}
Setting
$$\rho_i=(\sum_{j\in I} a^{ij})(4\pi \sum_l\gamma_l),\quad i\in I, $$
we can write (\ref{main-sys-5}) as
$$\Delta_g u_i^*+\sum_j a_{ij}\rho_j (\frac{h_j^*e^{u_j^*}}{\int_M h_j^*e^{u_j^*}}-1)=\sum_p 4\pi \gamma_l(\delta_{p_l}-1),\quad i\in I. $$
If $M$ is a torus ($\chi(M)=0$) and $\gamma_l\in \mathbb N^+$ we can compute the Leray-Schauder degree if $\sum_l \gamma_l$ is odd.
\begin{thm}\label{thm-spe}
Suppose $M$ is a torus, $\gamma_l\in \mathbb N^+$ and $\sum_{l}\gamma_l$ is odd. Then Leray-Schauder degree for (\ref{main-sys-5}) is $\frac 12 \Pi_{l=1}^N(1+\gamma_l)$.
\end{thm}

Here we would like to point out that if the topology of the manifold is trivial, Bartolucci \cite{bart-15} studied another delicate Liouville system and proved some existence results when the topological degree is zero. 

The main ideas of proofs in this article are motivated by a number of related works. One major difficulty comes from the ``partial blowup phenomenon", which means when a system is scaled according to the maximum of all its branches, some components disappear after taking the limit. One crucial step is to prove that no component is lost after scaling.  We call this a fully bubbling phenomenon. For this part we use the idea in \cite{linzhang2}. Another major difficulty comes from the non-simple blowup phenomenon. When a singular source happens to be a blowup point, it is possible to have a finite number of disjoint bubbling disks all tending to the singular source. Such a blowup picture is called ``non-simple blowup", studied by Kuo-Lin \cite{kuo-lin} and independently by Bartolucci-Tarantello \cite{bart-ta} for singular Liouville equations. In this article, using ideas in \cite{lwz-apde,lwzy-apde} we extend the results of Kuo-Lin, Bartolucci-Tarantello to Liouville systems and prove that the non-simple blowup phenomenon can only occur if the strength of the singular source is a multiple of $4\pi$.

Finally we would like to explain the role of $(H1)$ and $(H2)$ and how the blowup analysis of Liouville systems is different from that of Toda systems \cite{lwz-apde,lwzy-apde}. For Liouville systems, the total integration (energy) of global solutions belongs to a hypersurface \cite{lin-zhang-1}, which means the energy is not discrete. To rule out the difficulty caused by the abundance of energy we need to use $(H1)$ and $(H2)$ to prove that the profiles of bubbling solutions around different blowup points are the same. Moreover, there is almost no energy outside the bubbling disks. However, for Toda systems, even though the energy of global solutions is quantized, a major difficulty comes from the fact that there is a lot of energy outside bubbling disks. In \cite{lwzy-apde}, tools in algebraic geometry are used to prove that energy outside bubbling disks is also quantized.

The organization of this article is as follows: In section two we analyse the asymptotic behavior of solutions near a blowup point and we prove, using ideas in \cite{linzhang2} that the energy of $u_i^k$ must satisfy certain rules around different blowup points. In this section we also establish certain estimates for non-simple blowup points. Then in section three we prove all the main theorems. In particular the proof of degree counting theorems is by reducing the systems to Liouville equation and use the previous results of Chen-Lin \cite{chenlin2,chenlin3}.

\section{Asymptotic behavior around a singular source}

Since the proof of all the main theorems boils down to detailed analysis of locally defined blow up solutions, in this section we consider a locally defined Liouville system
\begin{equation}\label{sys-loca-1}
\Delta u_i^k+\sum_{j=1}^n a_{ij} h_j^k e^{u_j^k}=4\pi \gamma \delta_0, \quad i\in I,\quad \mbox{ in }\quad B_{\delta}\subset \mathbb R^2
\end{equation}
where $h_1^k,...,h_n^k$ are positive smooth functions on $B_{\delta}$ (the ball centered at the origin with radius $\delta>0$) with uniform bounds:
\begin{equation}\label{h-loc}
0<c_1\le h_i^k\le c_2,\quad \|h_i^k\|_{C^1}\le c_3, \quad i\in I,
\end{equation}
for $c_1,c_2,c_3>0$ independent of $k$. Let $\gamma>-1$ is the strength of $\delta_0$, $A=(a_{ij})_{n\times n}$ satisfy $(H1)$, $(H2)$, and we assume the uniform bound on the integral of $h_i^ke^{u_i^k}$ and its oscillation on $\partial B_{\delta}$ (the boundary of $B_{\delta}$) :
\begin{equation}\label{e-unit}
\int_{B_{\delta}}h_i^ke^{u_i^k}\le C
\end{equation}
\begin{equation}\label{bry-osci}
\max_i\max_{x,y\in \partial B_{\delta}}|u_i^k(x)-u_i^k(y)|\le C,
\end{equation}
for some $C$ independent of $k$. Then in this section we consider the case that the origin is the only blowup point in $B_{\delta}$: let
\begin{equation}\label{t-uk}
\tilde u_i^k(x)=u_i^k(x)-2\gamma \log |x|,\quad i\in I:=\{1,...,n\}
\end{equation}
and write the equation for $\tilde u^k=(\tilde u_1^k,...,\tilde u_n^k)$ as
\begin{equation}\label{e-tuk}
\Delta \tilde u_i^k+\sum_j a_{ij} |x|^{2\gamma} h_j^k e^{\tilde u_j^k}=0, \quad \mbox{ in }\quad B_{\delta}.
\end{equation}
Then we assume that
\begin{equation}\label{mk-muk}
M_k=\max_i \max_{x\in B_{\delta}}\frac{\tilde u_i^k(x)}{\mu} \quad \mbox{ where }\quad \mu=1+\gamma,
\end{equation}
tends to infinity:
\begin{equation}\label{b-u}
M_k\to \infty\quad \mbox{ and given}\,\, \epsilon\in (0,\delta),
\max_i\max_{x\in B_{\delta}\setminus B_{\epsilon}}u_i^k\le C(\epsilon)
\end{equation}
for some $C(\epsilon)>0$ independent of $k$.

In this case the profile of blowup solutions is more intriguing than that around a regular point. There are two possibilities: either
\begin{equation}\label{sphe-h}
\max_i \max_{x\in B_{\delta}}u_i^k(x)+2\log |x|\le C
\end{equation}
or along a subsequence
\begin{equation}\label{n-s-h}
\max_i\max_{x\in B_{\delta}}u_i^k(x)+2\log |x|\to \infty.
\end{equation}
We call the blowup phenomenon ``simple" if (\ref{sphe-h}) holds. Otherwise, if (\ref{n-s-h}) holds we use ``non-simple-blowup" to describe $u^k$.

\medskip

\subsection{Simple-blowup}

First we consider the case when (\ref{sphe-h}) holds.
Let
$$\tilde v_i^k(y)=\tilde u_i^k(\epsilon_k y)+2\mu \log \epsilon_k,  \quad
\mbox{where} \quad \epsilon_k=e^{-\frac 12 M_k}.$$
Then it is easy to verify that
$\tilde v_i^k\le 0$
and
\begin{equation}\label{e-tvk}
\Delta \tilde v_i^k(y)+\sum_j a_{ij} h_j(\epsilon_k y) |y|^{2\gamma} e^{\tilde v_j^k(y)}=0, \quad |y|\le \delta\epsilon_k^{-1}.
\end{equation}
Then we prove
\begin{lem}\label{lower-b}
\begin{equation}\label{tvk0}
\max_i \tilde v_i^k(0)\ge -C.
\end{equation}
\end{lem}

\noindent{\bf Proof:} From (\ref{b-u}) we  see that there exists $y_k\in B(0,\delta \epsilon_k^{-1})$ such that $\max_i \tilde v_i^k(y_k)=0$ and
$|y_k|=o(1)\epsilon_k^{-1}$. Let $r_k=|y_k|$ and
$$z_i^k(y)=\tilde v_i^k(r_ky)+2\mu \log r_k-c_0, \quad |y|\le 2, \quad i\in I$$
where $c_0$ is chosen to make $z_i^k\le -1$ (see (\ref{sphe-h})).
We write the equation of $z_i^k$ as
$$\Delta z_i^k+\frac{\sum_ja_{ij}|y|^{2\gamma}h_j^ke^{z_j^k+c_0}}{z_i^k}z_i^k=0,\quad |y|\le 2. $$
Using $z_i^k\le -1$ we see that
$|\frac{\sum_j a_{ij}h_j^ke^{z_i^k+c_0}}{z_i^k}|$ is bounded. Standard Harnack inequality for linear equation gives
\begin{equation}\label{har-z}
\max_{\partial B_1}(-z_i^k)\le C \min_{\partial B_1}(-z_i^k), \quad i\in I.
\end{equation}
Thus $\max_i\min\tilde v_i^k\ge -C$ on $\partial B_{r_k}$. Then (\ref{tvk0}) follows easily from standard maximum principle.
Lemma \ref{lower-b} is established. $\Box$

\medskip

The proof of Lemma \ref{lower-b} also implies that at least one component of $\tilde v_i^k$ is bounded below over any compact subset of $\mathbb R^2$, which means these components converge to a global function along a subsequence.
Thus we use $I_1$ to be the indexes of converging components.
In other words, for indexes not in $I_1$, the corresponding components tend to minus infinity over any fixed compact subset of $\mathbb R^2$.

Let $\tilde v_i$ be the limit of $\tilde v_i^k$ and we use
$$\sigma_i=\frac{1}{2\pi }\int_{\mathbb R^2} |y|^{2\gamma} e^{\tilde v_i}, \quad i\in I_1 $$
to denote the energy of $\tilde v_k$ in $\mathbb R^2$. Here for convenience we assumed $h_i^k(0)=1$, but this assumption is not essential.
Traditional method can be used to prove
$$\tilde v_i(y)=-m_i \log |y|+O(1),\quad |y|>1,\quad i\in I_1, $$
where $m_i=\sum_{j=1}^na_{ij}\sigma_j$.
For $i\in I_1$ we have
\begin{equation}\label{mi-big}
m_i>2\mu,\quad \mu=1+\gamma,\quad i\in I_1.
\end{equation}
Let $\sigma_i^k$ denote the energy of $u_i^k$ in $B_{\delta}$:
$$\sigma_i^k=\frac 1{2\pi}\int_{B_{\delta}}h_i^k|x|^{2\gamma}e^{\tilde u_i^k},\quad i=1,...,n,$$
Then it is immediate to observe that
$$\lim_{k\to \infty}\sigma_i^k\ge\sigma_i,\quad i\in I_1. $$
 Corresponding to $\sigma_i^k$ we set $m_i^k$ to be
 $$m_i^k=\sum_{j=1}^na_{ij}\sigma_j^k. $$
Before we proceed we extend (\ref{mi-big}) to all $i\in I$:

\begin{lem}\label{mi-big-1}
\begin{equation}\label{small-3}
m_i=\sum_{j\in I}a_{ij}\sigma_j>2\mu,\quad i\in I\setminus I_1.
\end{equation}
\end{lem}

\noindent{\bf Proof:} First we invoke a result from \cite{linzhang2}: For $A$ satisfying $(H1)$ and $(H2)$, $a_{ij}>0$ if $i\neq j$.
We prove (\ref{small-3}) by contradiction. Suppose $m=\min\{m_i,\quad i\in I\}\le 2\mu$. Then we immediately observe two facts: first $m_i>m$ for all $i\in I_1$ because $m_i>2\mu$ for $i\in I_1$. Second $m>0$ because $\sigma_i=0$ if $i\not \in I_1$ and $a_{ij}>0$ if $i\neq j$. Let $J=\{i\in I;\quad m_i=m\}$. Clearly $J$ is not empty, $I_1\cap J=\emptyset$ and we use $J_1$ to denote $I\setminus \{I_1\cup J\}$. Moreover we use $\bar m=\min\{m_i;\quad i\in I_1\cup J_1\}$. Clearly $\bar m>m$.
For each $i\in J$, we have $\sigma_i=0$ since $i\not \in I_1$. Thus
$$0=\sigma_i=\sum_ja^{ij}m_j=\sum_{i\in J}ma^{ij}+\sum_{j\in J_1\cup I_1}a^{ij}m_j. $$
Using $m_j>\bar m$ for $i\not \in J$ and $a_{ij}>0$ for $i\neq j$, we have
\begin{equation}\label{add-1}
0\ge m\sum_{j\in J}a^{ij}+\sum_{j\not\in J}a_{ij}\bar m=m\sum_{j\in I}a^{ij}+\sum_{j\not \in J}a^{ij}(\bar m-m).
\end{equation}
In view of $(H2)$, which includes $\sum_ja^{ij}\ge 0$, we see that equality in (\ref{add-1}) holds and
$$a^{ij}=0,\quad \forall i\in J \mbox{ and } \quad \forall j\in I\setminus J. $$
Thus $A^{-1}$ can be written as a block-diagonal form, which means $A$ can also be written as a block diagonal form (after possible rearrangement of indexes), which is a
contradiction to the irreducibility of $A$. (\ref{small-3}) and Lemma \ref{mi-big-1} are established. $\Box$

\medskip

The following lemma gives an estimate of the behavior of $u_i^k$ near $\partial B_{\delta}$:

\begin{lem}\label{local-est-1}
Let $M_k$ be defined in (\ref{mk-muk}) and $0$ be a simple blowup point of $u^k$, then we have
$$\sigma_i^k=\sigma_i+o(1),\quad i=1,2...,n, $$
$$\tilde u_i^k(x)=-m_i^k\log |x|-\frac{m_i^k-2\mu}{2}M_k+O(1),\quad x\in \partial B_{\delta}, \quad i\in I_1,$$
and
$$\tilde u_i^k(x)=-m_i^k\log |x|-\frac{m_i^k-2\mu}{2}M_k+(\tilde u_i^k(0)-\mu M_k)+O(1),\quad x\in \partial B_{\delta}, \quad i\not \in I_1.$$
\end{lem}

\begin{rem} Note that we use $o(1)$ to denote a quantity tending to $0$ as $k\to \infty$, and $O(1)$ to denote a quantity whose absolute value does not tend to infinity as $k\to \infty$. For $i\not \in I_1$, $\tilde u_i^k(0)-\mu M_k\to -\infty$. Also even for $i\not \in I_1$, $\lim_{k\to \infty}m_i^k>0$ because
$a_{ij}>0$ for $i\neq j$.
\end{rem}

\noindent{\bf Proof of Lemma \ref{local-est-1}:} As mentioned before, at least one component of $\tilde v_k$ converges uniformly over any fixed compact subset of $\mathbb R^2$.
Then it is easy to find $R_k\to \infty$ to make the following hold:
$$\frac{1}{2\pi}\int_{B_{R_k}}|y|^{2\gamma}h_i^k(\epsilon_ky)e^{\tilde v_i^k(y)}dy= \sigma_i+o(1),\quad i\in I_1, $$
$$\frac{1}{2\pi}\int_{B_{R_k}}|y|^{2\gamma}h_i^k(\epsilon_ky)e^{\tilde v_i^k(y)}dy=o(1),\quad i\in I\setminus I_1. $$

Let $\overline{\tilde v_i^k}(r)$ be the spherical average of $\tilde v_i^k$ on $\partial B_r$,
the differentiation of $\overline{\tilde v_i^k}(r)$ gives
$$\frac{d}{dr}\overline{\tilde v_i^k}(r)=\frac{1}{2\pi r}\int_{B_r}\Delta \tilde v_i^k=-\frac{1}{2\pi r}\int_{B_r}\sum_j a_{ij}h_j^k(\epsilon_ky)|y|^{2\gamma}e^{\tilde v_j^k}. $$
Since $a_{ij}\ge 0$ and all $m_i>2\mu$, it is easy to use Green's representation of $\tilde v_i^k$ to prove
$$\overline{\tilde v_i^k}(r)=\tilde v_i^k(y)+O(1),\quad |y|=r,\quad R_k\le r\le \frac{\delta}2\epsilon_k^{-1} $$
and
\begin{equation}\label{small-e}
\tilde v_i^k(y)\le -\overline{\tilde v_i^k}(R_k)-(2\mu+\delta_1)\log |y|+O(1),\quad |y|\ge R_k,\quad i\in I
\end{equation}
for some $\delta_1>0$ independent of $k$. Even though $\delta_1>0$ may be small, it leads to the smallness of the energy of $\tilde v_i^k$:
$$\int_{B_{\epsilon_k^{-1}\delta}\setminus B_r}|y|^{2\gamma} h_i^k(\epsilon_ky)e^{\tilde v_i^k}=O(r^{-\delta_1}) $$
Thus we can give an accurate estimate of the energy of $\tilde v_i^k$ as:
\begin{equation}\label{small-2}
\frac{1}{2\pi}\int_{B_r}|y|^{2\gamma}h_i^k(\epsilon_ky)e^{\tilde v_i^k}=m_i^k-O(r^{-\delta_1}),\quad R_k\le r\le \delta \epsilon_k^{-1}.
\end{equation}
By the smallness of the error term in (\ref{small-2}) and standard estimates from the Green's representation for $\tilde v_i^k$, we easily obtain
$$\tilde v_i^k(y)=-m_i^k\log |y|+O(1),\quad 1<|y|<\epsilon_k^{-1}\delta,\quad \forall i\in I_1. $$
The estimate for $\tilde v_i^k$ near infinity can be translated into the following estimate for $\tilde u_i^k$:
\begin{align}\label{u-ik-long}
\tilde u_i^k(x)&=\tilde v_i^k(y)-2\mu \log \epsilon_k \quad \mbox{ for } |x|=\delta, \quad |y|=\epsilon_k^{-1}\delta, \nonumber \\
=&-m_i^k \log |y|-2\mu (-\frac 12 M_k)+O(1) \nonumber \\
=&-m_i^k \log |x|+m_i^k \log \epsilon_k +\mu M_k+O(1) \nonumber \\
=&-m_i^k \log |x|-\frac{m_i^k-2\mu}{2}M_k +O(1).
\end{align}
Thus the estimate for $i\in I_1$ for $u_i^k$ is established.

It is also straight forward to prove that for all $i\not \in I_1$,
$$\int_{B(0,\epsilon_k^{-1}\delta)\setminus B_r}|y|^{2\mu}h_i^k(\epsilon_k y)e^{\tilde v_i^k(y)}=O(r^{-\delta_1}),\quad R_k\le r\le \epsilon_k^{-1}\delta. $$
With this estimate the behavior of $\tilde v_i^k$ for $i\not \in I_1$ can be written as
$$\tilde v_i^k(y)=\tilde v_i^k(0)-(\sum_{j\in I_1} a_{ij}\sigma_j^k)\log |y|+O(1), \quad i\not \in I_1. $$
Consequently for $\tilde u_i^k$ we have, for $|x|=\delta$ and $|y|=\epsilon_k^{-1}|x|$,
\begin{align*}
\tilde u_i^k(x)&=\tilde v_i^k(y)-2\mu \log \epsilon_k\\
&=\tilde v_i^k(0)-(\sum_{j\in I_1} a_{ij} \sigma_j^k) \log |y| +\mu M_k, \\
&=-m_i^k\log |x|-\frac{m_i^k-2\mu}{2}M_k+(\tilde u_i^k(0)-\mu M_k).
\end{align*}
Lemma \ref{local-est-1} is established. $\Box$

\begin{rem}\label{global-pi}
Even though some components of $(\sigma_1,....,\sigma_n)$ may be zero, $(\sigma_1,...,\sigma_n)$ still satisfies the standard Pohozaev identity:
\begin{equation}\label{pi-sig}
\sum_{i,j\in I} a_{ij}\sigma_i\sigma_j=4\sum_i (1+\gamma)\sigma_i.
\end{equation}
The derivation of (\ref{pi-sig}) is standard and we mention the argument here for the convenience of readers. The Pohozaev identity for $u^k$ on $\Omega$ is
\begin{eqnarray*}
&\sum_{i\in I}(\int_{\Omega}(x\cdot \nabla h_i^k)e^{u_i^k}+2h_i^ke^{u_i^k})\\
&=\int_{\partial \Omega}\bigg (\sum_i(x\cdot \nu)h_i^ke^{u_i^k}+\sum_{i,j}a^{ij}(\partial_{\nu}u_j^k(x\cdot \nabla u_i^k)-\frac 12(x\cdot \nu)(\nabla u_i^k\cdot \nabla u_j^k))\bigg )
\end{eqnarray*}
Setting $\Omega=B_{\delta}\setminus B_{\epsilon}$ and let $\epsilon\to 0$, we have
\begin{eqnarray*}
\delta \int_{\partial B_{\delta}}\sum_{ij\in I}a^{ij}(\partial_{\nu}u_i^k\partial_{\nu}u_j^k-\frac 12\nabla u_i^k\cdot \nabla u_j^k)
+\sum_{i\in I}\delta \int_{\partial B_{\delta}}h_i^ke^{u_i^k}\\
=2\sum_{i\in I}\int_{B_{\delta}}h_i^ke^{u_i^k}+\sum_{i\in I}\int_{B_{\delta}}(x\cdot \nabla h_i^k)e^{u_i^k}+4\pi \sum_{ij\in I}a^{ij}\gamma^2.
\end{eqnarray*}
where we have used
$$\nabla u_i^k=2\gamma x/|x|^2+ \mbox{ a bounded function } $$
near the origin. In order to evaluate other terms we can use standard elliptic estimate to obtain
$$\nabla u_i^k(x)=(\sum_{i\in I}a_{ij}\sigma_j-2\gamma+o(1))/|x|,\quad |x|=\delta. $$
Then (\ref{pi-sig}) follows from direct computation. We refer the readers to \cite{lwz-apde} and \cite{lin-zhang-1} for more detailed computation.
\end{rem}

\begin{rem} If the blowup point $p$ is not a singular source, the scaling is centered at $p_k\to p$ where the maximum of $\tilde u_i^k$ is attained. In this case we have $\max_iv_i^k(0)=0$ and the non-simple blow-up does not happen.
\end{rem}

\medskip

\subsection{The comparison of blowup solutions around different blowup points}

Under the same context as in the previous subsection, we establish the
 following lemma which compares the behavior of solutions outside bubbling disks.

\begin{lem}\label{two-bubbles}
Let $p$ and $q$ be two disjoint blowup points of
$$\Delta u_i^k+\sum_j a_{ij}h_j^ke^{u_j^k}=4\pi \gamma_p \delta_p+4\pi \gamma_q \delta_q, \quad i\in I $$
in $\Omega\subset \subset\mathbb R^2$ where $p,q\in \Omega$, $\gamma_p,\gamma_q>-1$. Suppose the assumption on $h_i^k$ is the same as before: Uniformly bounded above and below by positive constants and uniformly bounded in $C^1$ norm. And we also have the uniform bound energy and finite oscillation assumptions:
$$\int_{\Omega}h_i^ke^{u_i^k}\le C, \quad \max_{x,y\in \partial \Omega} |u_i^k(x)-u_i^k(y)|\le C. $$
 We use $(\sigma_1^k,...,\sigma_n^k)$ and $(\bar \sigma_1^k,...,\bar \sigma_n^k)$ to denote the integration of $u^k$ in $B(p,\delta)$ and $B(q,\delta)$, respectively:
$$\sigma_i^k=\frac 1{2\pi}\int_{B(p,\delta)}h_i^ke^{u_i^k},\quad \bar \sigma_i^k=\frac 1{2\pi}\int_{B(q,\delta)}h_i^ke^{u_i^k}. $$
If $p$ or $q$ is a regular point instead of a singular source, we have $\gamma_p=0$ or $\gamma_q=0$. Correspondingly we set
$$m_i^k=\sum_{j\in I}a_{ij}\sigma_j^k,\quad \bar m_i^k=\sum_{j=1}^n a_{ij}\bar \sigma_j^k. $$
Assume in addition that
 $$u_i^k|_{\partial B(p,\delta)}=u_i^k|_{\partial B(q,\delta)}+O(1). $$
 Then if $p$ and $q$ are both simple blowup points, we have
 \begin{equation}\label{same-1}
 \frac{\mu_q}{\mu_p}\lim_{k\to \infty}\sigma_i^k=\lim_{k\to \infty}\bar \sigma_i^k,\quad i\in I.
 \end{equation}
\end{lem}

\begin{rem} If $p$ or $q$ is a regular point, it is a simple blowup already.
\end{rem}

\noindent{\bf Proof of Lemma \ref{two-bubbles}:}

Since $p$ and $q$ can be a singular source or a regular point on the manifold, we use $\mu_p=1+\gamma_p$ if $p$ is a singular source. Otherwise $\mu_p=1$.
Let $M_k=\max_{i\in I}\tilde u_i^k(x)/\mu_p$ for $x\in B(p,\delta)$ and $\bar M_k=\max_{i}\tilde u_i^k(x)/\mu_q$ in $B(q,\delta)$, where $\tilde u_i^k$ is $u_i^k$ minus a corresponding logarithmic term in local coordinates.  Suppose $M_k$ is attained at $p_k$ that tends to $p$ and $\bar M_k$ is attained at $q_k$ that tends to $q$. Using Lemma \ref{local-est-1} we have, for $i\in I$
\begin{equation}\label{comp-1}
\frac{m_i^k-2\mu_p}{2}M_k+(\mu_p M_k-\tilde u_i^k(p_k))=\frac{\bar m_i^k-2\mu_q}{2}\bar M_k+(\mu_q \bar M_k-\tilde u_i^k(q_k))+O(1).
\end{equation}
Here we further remark that, say around $p$, if the first $l$ components of $\tilde u^k$ converge to a system of $l$ equations after scaling, $\mu_p M_k-\tilde u_i^k(p_k)$ are uniformly bounded for $1\le i\le l$. In this case $\mu_p M_k-\tilde u_i^k(p_k)$ can be replaced by $O(1)$. For $i>l$, $\mu_p M_k-\tilde u_i^k(p_k)$ tends to infinity. The right hand side of (\ref{comp-1}) can be understood similarly. For each $i\in I$, if
$$\mu_p M_k-\tilde u_i^k(p_k)>\mu_q\bar M_k-\tilde u_i^k(q_k), $$
we let
$$l_i^k=(\mu_pM_k-\tilde u_i^k(p_k))-(\mu_q\bar M_k-\tilde u_i^k(q_k)),\quad \mbox{and}\quad \bar l_i^k=0. $$
On the other hand, if
$$\mu_pM_k-\tilde u_i^k(p_k)\le \mu_q\bar M_k-\tilde u_i^k(q_k), $$
we let
$$l_i^k=0,\quad \mbox{and}\quad \bar l_i^k=(\mu_q\bar M_k-\tilde u_i^k(q_k))-(\mu_p M_k-\tilde u_i^k(p_k)). $$

Set
$$I_1:=\{i\in I;\quad \lim_{k\to \infty}\frac{l_i^k}{M_k}>0\}\quad, \mbox{ and } \quad I_2:=\{i\in I;\quad \lim_{k\to \infty} \frac{\bar l_i^k}{M_k}>0\}. $$
 It is easy to observe that $I_1\cap I_2=\emptyset$. We claim that $I_1=\emptyset$, which is now proved by contradiction:

Suppose $I_1\neq \emptyset$, then we consider two cases: $I_2\neq \emptyset$ and  $I_2=\emptyset$.

\medskip

\noindent{\bf Case one: $I_2\neq \emptyset$}.

Let
$$\lambda=\lim_{k\to \infty}\frac{M_k}{\bar M_k}, \quad \delta_i=\lim_{k\to \infty} \frac{l_i^k}{\bar M_k},\quad \bar \delta_i=\lim_{k\to \infty}\frac{\bar l_i^k}{\bar M_k}. $$
We claim that these limits exist along a subsequence. Indeed, using the definition of $l_i^k$ and $\bar l_i^k$ (\ref{comp-1}) can be written as
$$\frac{m_i^k-2\mu_p}{2}\frac{M_k}{\bar M_k}+\frac{l_i^k}{\bar M_k}=\frac{\bar m_i^k-2\mu_q}{2}+\frac{\bar l_i^k}{\bar M_k}+o(1). $$
Take $i\in I_1$, the right hand side tends to $\frac{\bar m_i-2\mu_q}{2}$, which means along a subsequence, the two terms on the left hand side tend to $\frac{m_i-2\mu_p}2$ and $\delta_i$, respectively( we use $\sigma_i$ to denote the limit of $\sigma_i^k$. $m_i$,$\bar m_i$,$\bar \sigma_i$ are understood in a similar fashion). On the other hand for $j\in I_2$, the left hand side tends to $\frac{m_i-2\mu_p}{2}\lambda$, which forces the right hand side to converge to $\frac{\bar m_i-2\mu_q}{2}+\bar \delta_i$ along a subsequence. Now (\ref{comp-1}) leads to
\begin{equation}\label{comp-2}
\lambda\frac{m_i-2\mu_p}{2}+\delta_i=\frac{\bar m_i-2\mu_q}{2}+\bar \delta_i,\quad \forall i\in I.
\end{equation}

Here we recall that  $\delta_i>0$ in $I_1$ and $\bar \delta_i>0$ in $I_2$. We also will use $\sigma_i\delta_i=0$ for all $i$.
From $\bar \sigma_i=0$ in $I_2$, we have
$$0=\bar \sigma_i=\sum_{j\in I_2}a^{ij}\bar m_j+\sum_{j\not \in I_2}a^{ij}\bar m_j. $$
Since $A$ is irreducible, there exist $i\in I_2$ and $j\not \in I_2$ such that $a^{ij}>0$.
Multiplying $\bar \delta_i$ on both sides and taking the summation for $i\in I_2$, we have
$$\sum_{i,j\in I_2}a^{ij}\bar m_i\bar \delta_j<0. $$
So trivially there exists $\tilde i\in I_2$ such that
\begin{equation}\label{add-2}
\sum_{j\in I_2}a^{\tilde ij}\bar \delta_j<0.
\end{equation}
From the comparison of the $\tilde i$th component, we have
$$\lambda\sum_j a^{\tilde ij}(\frac{m_j-2\mu_p}2)+\sum_j a^{\tilde i j}\delta_j=\frac{\bar \sigma_{\tilde i}}2-\sum_j a^{\tilde ij}\mu_q+\sum_j a^{\tilde i j} \bar \delta_j. $$
The second term on the left is nonnegative because $\delta_i=0$ if $i\in I_2$ and $a^{\tilde ij}\ge 0$ if $\tilde i\neq j$. The first term on the right is $0$, the last term on the right is negative. Thus
the equation above is reduced to
$$\frac{\lambda}2\sigma_{\tilde i}-\lambda \mu_p\sum_j a^{\tilde ij}<-\sum_{j}a^{\tilde ij}\mu_q. $$
Since $\sigma_{\tilde i}\ge 0$, the strict inequality and $(H2)$ imply $\sum_ja^{\tilde ij}>0$, thus we have
$$\lambda>\mu_q/\mu_p. $$
On the other hand the same argument applied to $i\in I_1$ gives
$$\lambda<\frac{\mu_q}{\mu_p}. $$
 Thus  this case ($I_1\neq \emptyset$, $I_2\neq \emptyset$ ) is ruled out.

\medskip

Next under the assumption $I_1\neq \emptyset$ we consider the case that $I_2=\emptyset$.

Since all $\bar \delta_i=0$ we have
$$\frac{m_i-2\mu_p}{2}\lambda+\delta_i=\frac{\bar m_i-2\mu_q}2,\quad i\in I. $$
Using this expression in
$$\sum_{ij}a^{ij}(\frac{\bar m_i-2\mu_q}{2})(\frac{\bar m_j-2\mu_q}2)=\sum_{ij}a^{ij}\mu_q^2, $$
which is equivalent to the Pohozaev identity for $(\bar \sigma_1,...,\bar \sigma_n)$ (see Remark \ref{global-pi})
we have
\begin{equation}\label{i1i2-1}
\lambda^2\sum_{ij}a^{ij}\mu_p^2+2\lambda\sum_{ij}a^{ij}(\frac{m_i-2\mu_p}2)\delta_j+\sum_{i,j\in I_1}a^{ij}\delta_i\delta_j=\sum_{ij}a^{ij}\mu_q^2
\end{equation}
where we have used
$$\sum_{ij}a^{ij}(\frac{m_i-2\mu_p}2)(\frac{m_j-2\mu_p}2)=\sum_{ij}a^{ij}\mu_p^2. $$

The second term on the left hand side of (\ref{i1i2-1}) can be written as
$$\lambda( \sum_j\sigma_j\delta_j-2\mu_p\sum_j(\sum_i a^{ij})\delta_j), $$
which is nonpositive because $\sigma_i\delta_i=0$ and $\sum_ia^{ij}\ge 0$. We further claim that the third term on the left hand side of (\ref{i1i2-1}) is nonpositive. This is because all the eigenvalues of $(a^{ij})_{I_1\times I_1}$ are non-positive. This is  proved in \cite{linzhang2} and we include it here for convenience: Without loss of generality we assume $I_1=\{1,...,i_0\}$ and let $\mathbb{F}=(a^{ij}){i_0\times i_0}$ for $i,j\in I_1$. Let $\mu$ be the largest eigenvalue of $\mathbb{F}$ and $\eta=(\eta_1,...,\eta_{i_0})$ be an eigenvector corresponding to $\mu$. Here $\eta$ is the vector that attains
$$\max_{v\in \mathbb R^{i_0}}v^T\mathbb{F}v,\quad v^Tv=1. $$
Since $a^{ij}\ge 0$ for all $i\neq j$, we can choose $\eta_i\ge 0$ for all $i\in I_1$. For each $i\in I_1$,
$$0=\sigma_i=\sum_{j\in I_1}a^{ij}m_j+\sum_{j\not \in I_1}a^{ij}m_j. $$
Thus by (H2)
$$\sum_{j\in I_1}a^{ij}m_j\le 0,\quad i\in I_1. $$
Multiplying both sides by $\eta_i$ and taking summation on $i$, we have
$$0\ge \sum_{i,j\in I_1}a^{ij}\eta_im_j=\sum_{j\in I_1}\mu \eta_jm_j. $$
Using $\eta_i\ge 0$ (at least one of them is strictly positive) and $m_i>0$ for $i\in I_1$, we have $\mu\le 0$.

\medskip

 Thus from (\ref{i1i2-1}) we have
$$\lambda\ge \frac{\mu_q}{\mu_p}. $$
Note that we have used $\sum_{ij}a^{ij}>0$ because otherwise $A^{-1}$ would not be invertible.
Next using the proof of (\ref{add-2}) we can find some $i\in I_1$ such that $\sum_{j\in I_1}a^{ij}\delta_j<0$. For this $i$, from
$$\sum_{j}a^{ij}(\frac{m_j-2\mu_p}2\lambda+\delta_j)=\sum_ja^{ij}\frac{\bar m_j-2\mu_q}2 $$
we write it as
$$-\sum_ja^{ij}\mu_p\lambda+\sum_{j\in I_1}a^{ij}\delta_j=\bar \sigma_i/2-\sum_j a^{ij}\mu_q, $$
where we have use $\sigma_i=0$ for $i\in I_1$.
Using $\bar \sigma_i\ge 0$ and $\sum_{j\in I_1}a^{ij}\delta_j<0$ we have
$$\lambda< \frac{\mu_q}{\mu_p}. $$
Therefore  this case ($I_1\neq \emptyset, I_2=\emptyset$) is also ruled out. We have proved that $I_1=\emptyset$. In a similar manner $I_2=\emptyset$ can also be established.

Finally using
\begin{equation}\label{ratio-2}
\lambda\frac{m_i-2\mu_p}2=\frac{\bar m_i-2\mu_q}2,\quad i\in I,
\end{equation}
in the Pohozaev identity for $(\bar \sigma_1,...,\bar \sigma_n)$ we have
\begin{equation}\label{ratio-1}
\lambda=\lim_{k\to \infty}\frac{M_k}{\bar M_k}=\frac{\mu_q}{\mu_p}.
\end{equation}
Using (\ref{ratio-1}) in (\ref{ratio-2}) we further have
\begin{equation}\label{ratio-3}
\frac{\mu_q}{\mu_p}m_i=\bar m_i,\quad \frac{\mu_q}{\mu_p}\sigma_i=\bar \sigma_i, \quad i\in I.
\end{equation}

Lemma \ref{two-bubbles} is established $\Box$
\medskip

Finally we deduce the asymptotic behavior of $u^k$ when non-simple-blowup occurs.

\medskip

\subsection{Non-simple blowup}
Now we consider the second possibility, the non-simple blowup. This phenomenon happens when (\ref{n-s-h}) holds. Recall that $u^k=(u_1^k,...,u_n^k)$ satisfies
(\ref{sys-loca-1}). If (\ref{n-s-h}) holds, a standard selection process ( \cite{lwz-apde} ) determines a finite number of bubbling disks: $B(p_l^k,r_l^k)$ for $l=1,...,N$ where $p_l^k$ are local maximums of some $u_i^k$ and $r_l^k$s are determined as follows: Scale $u^k$ with respect to the maximum of $\max_iu_i^k(p_l^k)$, then the system converges to a possibly smaller global system with finite energy. Note that we use $B(p,\delta)$ to denote the ball centered at $p$ with radius $\delta$. Then it is easy to choose $R_k\to \infty$ such that the integral of the scaled functions over $B(0,R_k)$ is only $o(1)$ different from the energy of entire solutions. Scaling back to $u_k$ we have that the integral of $e^{u_i^k}$ over $B(p_l^k,r_l^k)$ is $o(1)$ different from the energy of its global limit. Moreover, if we use $(\sigma_{l1}^k,...,\sigma_{ln}^k)$ to denote the energy in
$B(p_l^k,r_l^k)$ we have
$$\sum_{i,j\in I}a_{ij}\sigma_{li}^k\sigma_{lj}^k=4\sum_i\sigma_{li}^k+o(1). $$

Here we shall invoke some argument in \cite{lwz-apde}. The main result in this part is:

\begin{prop}\label{nu-int}
If (\ref{n-s-h}) holds, $\mu\in \mathbb N^+$.
\end{prop}

First we mention the following simple lemma:

\begin{lem}\label{one-bubble} Let $A=(a_{ij})_{n\times n}$ be a matrix that satisfies (H1).
Suppose $(\sigma_1^{(1)},....,\sigma_n^{(1)})$ and $(\sigma_1^{(2)},....,\sigma_n^{(2)})$ are two vectors with nonnegative components. If they
both satisfy
$$\sum_{i,j}a_{ij}\sigma_i^{(l)}\sigma_j^{(l)}=4\mu\sum_{i=1}^n\sigma_i^{(l)} $$
for $l=1,2$ and some $\mu>0$. Then if
\begin{equation}\label{poho-2e}
\sum_{j=1}^n a_{ij}\sigma_j^{(1)}>2\mu, \quad i=1,...,n
\end{equation}
and
$$\sigma_i^{(2)}\ge \sigma_i^{(1)}\quad i=1,...,n. $$
Then
$$\sigma_i^{(1)}=\sigma_i^{(2)},\quad i=1,...,n. $$
\end{lem}

\noindent{\bf Proof of Lemma \ref{one-bubble}:}  The proof is immediate. Let $s_i=\sigma_i^{(2)}-\sigma_i^{(1)}$. Then $s_i\ge 0$. The difference between the two equations in (\ref{poho-2e})
gives
$$\sum_{i,j}(a_{ij}\sigma_j^{(1)}-2\mu)s_i+\sum_{i,j}a_{ij}s_is_j=0. $$
By the assumption $(H1)$ and the nonnegativity of $s_i$ we have $s_i=0$ for all $i$. $\Box$

\medskip

\noindent{\bf Proof of Proposition \ref{nu-int}:} First we use

$$\Sigma_k=\{0,p_1^k,...,p_N^k\} $$
to denote the set of blowup points and the origin. Note that there may also be a bubbling disk centered at the origin, as described in Lemma \ref{lower-b}.
Here we invoke the definition of group in \cite{lwz-apde}. If a few bubbling disks are of comparable distance to one another and are much further to other bubbling disks, the set of these bubbling disks ( that of comparable distance to one another) is called a group. See \cite{lwzy-apde,lwz-apde} for more detailed discussions.
For example, $p_1^k,p_2^k,p_3^k$ are called in a group if $dist(p_1^k,p_2^k)\sim dist(p_1^k,p_3^k)\sim dist(p_2^k,p_3^k)$ and
$$dist(p_1^k,q)/dist(p_1^k,p_2^k)\to \infty $$
for any $q\in \Sigma_k\setminus \{p_1^k,p_2^k,p_3^k\}$.

Now we  make two important observations: First, there is no group far away from the origin. The reason is if there were
such a group, say $B(p_k,l_k)$ and $B(q_k,l_k)$ belong to a group and $dist(0,p_k)/dist(p_k,q_k)\to \infty$. First by the argument of Lemma \ref{mi-big} and Lemma \ref{local-est-1} all the components of $u_i^k$ have faster decay than harmonic function near $\partial B(p_k,l_k)$ and $\partial B(q_k,l_k)$: in precise terms, if we use $(\sigma_{p1}^k,...,\sigma_{qn}^k)$ and $(\sigma_{q1}^k,...,\sigma_{qn}^k)$ to denote the energy in $B(p_k,l_k)$ and $B(q_k,l_k)$, respectively, we have
$$\sum_{ij}a_{ij}\sigma_{pi}^k\sigma_{pj}^k=4\sum_i\sigma_{pi}+o(1), $$
and
$$\sum_{ij}a_{ij}\sigma_{qi}^k\sigma_{qj}^k=4\sum_i\sigma_{qi}+o(1). $$
Moreover, as in Lemma \ref{mi-big}
\begin{equation}\label{fast-p}
m_{pi}^k:=\sum_{ij}a_{ij}\sigma_{pj}^k>2,\quad m_{qi}^k:=\sum_ja_{ij}\sigma_{qj}^k>2,\quad \forall i\in I.
\end{equation}
Let $d_k$ be the distance from $p_k$ to the nearest member in $\Sigma_k$ not in the group of $p_k$ and $q_k$. Then (\ref{fast-p}) means all components of $u^k$ decay so fast that there is little energy in $B(p_k,d_k/2)\setminus (B(p_k,l_k)\cup B(q_k,l_k))$.
 Looking at the average of $u_i^k$ it is easy to find $\bar l_k\le d_k/2$
which satisfies
$$\bar l_k/l_k\to \infty, \bar l_k=o(1)dist(p_k, \Sigma_k\setminus \mbox{ the group of } p_k). $$
And on $\partial B(p_k,\bar l_k)$ we still have
\begin{equation}\label{pi-c}
u_i^k(x)+2\log l_k\to -\infty, \quad i\in I.
\end{equation}
From (\ref{pi-c}) it is easy to use the Green's representation formula to evaluate the Pohozaev identity and obtain (see \cite{lwz-apde})
\begin{equation}\label{add-pi}
\sum_{i,j}a_{ij} \sigma_{\bar li}^k\sigma_{\bar lj}^k=4\sum_{i}\sigma_{\bar li}^k+o(1), 
\end{equation}
where $\sigma_{\bar li}^k=\frac{1}{2\pi}\int_{B(p_k,\bar l_k)}h_i^ke^{u_i^k}$. Since $(\sigma_{\bar l1}^k,...,\sigma_{\bar ln}^k)$ and
$(\sigma_{p1}^k,...,\sigma_{pn}^k)$ satisfy the same equation but $\sigma_{\bar li}^k\ge \sigma_{pi}^k+\sigma_{qi}^k$, by Lemma \ref{one-bubble} we easily get a contradiction. Here we briefly review how (\ref{pi-c}) leads to (\ref{add-pi}). Roughly speaking (\ref{pi-c}) means the value of $u_i^k$ is very small on $\partial B(p_k,\bar l_k)$ and by Harnack inequality, most energy of $u_i^k$ in $B(p_k,\bar l_k)$ is concentrated near $p_k$, which implies that all the derivatives of $u_i^k$ is very easy to estimate on $\partial B(p_k,\bar l_k)$. The evaluation of the derivatives of $u_i^k$ and the smallness of $e^{u_i^k}$ on $\partial B(p_k,\bar l_k)$ lead to (\ref{add-pi}).

\medskip

The second main observation is that for the group containing the orgin, there is no bubbling disk centered at the origin. In other words, if there is a group that contains the origin, it has to be case that there are finitely many bubbling disks, say $B(p_1^k,r_1^k)$,...,$B(p_l^k,r_l^k)$, with $p_1^k$,...,$p_l^k$ all of comparable distance to the origin and there is no bubbling disk centered at the origin. This fact is also proved by contradiction. Suppose around the origin there is a bubbling disk whose energy is $(\sigma_1^k,...\sigma_n^k)$. We have already known that
$$\sum_{ij}a_{ij}\sigma_i^k\sigma_j^k=4\mu\sum_i \sigma_i^k+o(1). $$
If there is another bubbling disk, say $B(p_1^k,l_k)$ in the group, we can find $\bar l_k$ such that $B(0,\bar l_k)$ encloses all the bubbling disks in this group and
$\bar l_k$ is less than half of the distance from $0$ to any member in $\Sigma_k$ outside the group. The fast decay property as before also gives
$$u_i^k(x)+2\log \bar l_k\to -\infty, \quad x\in \partial B(0,\bar l_k). $$
Using the same argument as in \cite{lwz-apde} we have
$$\sum_{ij}\sigma_{\bar li}^k\sigma_{\bar lj}^k=4\mu\sum_i\sigma_{\bar li}^k+o(1), $$
where
$\sigma_{\bar i}^k=\frac{1}{2\pi}\int_{B(0,\bar l_k)}h_i^ke^{u_i^k}$. Since $\sigma_{\bar li}$ is significantly greater than $\sigma_i^k$ for at least one component, Lemma \ref{one-bubble} gives a contradiction as before.

\medskip

 By the two observations before we only need to consider the case that there are
finitely many bubbling disks around the origin and their centers are of comparable distance to the origin. Suppose these local maximums are $p_1^k$,...,$p_N^k$, and we suppose $|p_{t}^k|\sim \delta_k$.

Let
$$\Lambda_k=\max_i\max_x  u_i^k(x)+2\log |x|. $$
Without loss of generality we suppose $\Lambda_k$ is attained at $p_{1,k}$. Let $\delta_k=|p_{1,k}|$ and
\begin{equation}\label{u-vik}
v_i^k(y)=u_i^k(p_{1,k}+\delta_ky)+2\log \delta_k, \quad i\in I.
\end{equation}
It is immediate to observe that the domain of $v_i^k$ contains $B(0, \delta \delta_k^{-1})$ for some small $\delta>0$. Standard selection process can be employed to obtain finite bubbling disks centered at $p_{2,k}$,...,$p_{N,k}$ such that not only $|p_{j,k}|\sim \delta_k$, but also $|p_{m,k}-p_{l,k}|\sim \delta_k$ for all $l\neq m$. Let $z_l^k$ be the images of $p_{l,k}$ by the scaling in (\ref{u-vik}). Then clearly $z_1^k$ is the origin and the distance between any two $z_l^k$s is comparable to $1$. So we assume, $B(z_l^k,\delta)$ are mutually disjoint for some small $\delta>0$.
The definition of $v_i^k$ clearly implies that
$$\max_i\max_{B_{\delta}}v_i^k=\Lambda_k. $$
Let $I_1$ be the set of convergent components after scaling according to the maximum of all components. Then using previous discussion we have
$$
\left\{\begin{array}{ll}
v_i^k(y)=-m_i^k\log |y|-\frac{m_i^k-2}2\Lambda_k+O(1),\quad i\in I_1, \quad y\in \partial B_{\delta},\\
v_i^k(y)=-m_i^k\log |y|-\frac{m_i^k-2}2\Lambda_k+v_i^k(0)-\Lambda_k+O(1),\quad i\not\in I_1, \quad y\in \partial B_{\delta}.
\end{array}
\right.
$$
for some $\delta>0$. Here we use $(\sigma_1^k,...,\sigma_n^k)$ and $(m_1^k,...,m_n^k)$ to denote the energy around $p_{1,k}$:
$$\sigma_i^k=\frac 1{2\pi}\int_{B(0,\delta)} h_i^k(p_{1,k}+\delta_ky)e^{ v_i^k}, \quad i\in I, \quad m_i^k=\sum_ja_{ij}\sigma_j^k. $$

If we use $(\bar \sigma_1^k,...,\bar \sigma_n^k)$ and $(\bar m_1^k,...,\bar m_n^k)$ to denote energy around another bubbling disk in this group. Lemma \ref{two-bubbles} gives
$$\lim_{k\to \infty} \sigma_i^k=\lim_{k\to \infty}\bar \sigma_i^k,\quad i\in I. $$
The Pohozaev identity for $(\sigma_1^k,...,\sigma_n^k)$ is
\begin{equation}\label{poh-9}
\sum_{ij}a_{ij}\sigma_i^k\sigma_j^k=4\sum_i \sigma_i^k+o(1).
\end{equation}
 The equation for $(\bar \sigma_1^k,...,\bar \sigma_n^k)$ is the same. If we use $\bar \Lambda_k$ to denote the maximum around the bubbling disk that $\bar \sigma_i^k$ represents, the proof of Lemma \ref{two-bubbles} gives
$$\Lambda_k/\bar \Lambda_k=1+o(1). $$

Let $\sigma_i=\lim_{k\to \infty}\sigma_i^k$. Then $(\sigma_1,...,\sigma_n)$ satisfies
$$\sum_{ij}a_{ij}\sigma_i\sigma_j=4\sum_i \sigma_i. $$
On a fast decay radius that encloses all bubbling disks in the group round the singular source, we have
$$\sum_{ij}a_{ij}(N\sigma_i)(N\sigma_j)=4\mu\sum_i (N\sigma_i). $$
Thus $\mu=N$ ( that is $\gamma=N-1$) and Proposition \ref{nu-int} is established. $\Box$

\medskip

Next we derive the asymptotic behavior of $u_i^k$ on $\partial B_{\delta}$ for some $\delta>0$ small if the non-simple blowup phenomenon occurs. Recall that $\delta_k$ is the distance from $0$ to a local maximum of $u_i^k$. Here we abuse the notation of $v_i^k$ by defining it slightly differently:
$$v_i^k(y)=u_i^k(\delta_ky)+2\log \delta_k,\quad i\in I. $$
Then we have
$$
\Delta v_i^k(y)+\sum_{j}a_{ij}h_j^k(\delta_ky)e^{v_j^k}=4\pi \gamma \delta_0, \quad |y|<\delta \delta_k^{-1}, \quad i\in I.
$$
If we use $\bar v_i^k(r)$ to denote the spherical average of $v_i^k$ at $\partial B_r$, we have, for $r>>1$ ( so $B_r$ contains all the $N$ bubbling disks around the origin),
$$\frac{d}{dr}\bar v_i^k(r)=-\frac{1}{r}(\frac{1}{2\pi}\int_{B_r}a_{ij} h_j^ke^{v_j^k}-2\gamma)$$

Thus based on the asymptotic behavior of $v_i^k$ around each of the $N$ bubbling disks, we have
$$\frac{d}{dr}\bar v_i^k(r)=\frac{-N m_i^k+2\gamma+o(r^{-\delta_1})}{r} $$ for some $\delta_1>0$.
So for $r\sim \delta_k^{-1}$ we have, for $i\in I_1$,
$$v_i^k(y)=-\frac{m_i^k-2}2\Lambda_k+(-N m_i^k+2\gamma)\log \delta_k^{-1}+O(1),\quad |y|\sim \delta_k^{-1}. $$
Using $\gamma=N-1$ and the definition of $v_i^k$ in (\ref{u-vik}), we have
\begin{eqnarray}\label{a-i1}
&u_i^k|_{\partial B(p,\delta)}
=v_i^k|_{\partial B(0,\delta \delta_k^{-1})}+2\log \delta_k  \\
&=-\frac{m_i^k-2}2\Lambda_k-(\frac{m_i^k-2}2)2N\log \delta_k^{-1}+O(1), \nonumber\\
&=-\frac{m_i^k-2}2(\Lambda_k+2N\log \delta_k)+O(1),\quad i\in I_1. \nonumber
\end{eqnarray}

For $i\in I\setminus I_1$, we have
\begin{equation}\label{a-i2}
u_i^k|_{\partial B(p,\delta)}=-\frac{m_i^k-2}2(\Lambda_k+2N\log \delta_k)-N_k,
\end{equation}
for some $N_k=\Lambda_k-v_i^k(0)+O(1)\to \infty$.

\medskip

From (\ref{a-i1}) and (\ref{a-i2}) we see that even if the non-simple blowup phenomenon happens around a singular source, still the argument of Lemma \ref{two-bubbles} can be applied to compare the energy of two blowup points, regardless of they are simple or not. Thus under the same context of Lemma \ref{two-bubbles} except that we remove the simple-blowup requirement, we still have
\begin{equation}\label{same-e-2}
\frac{\sigma_{pi}}{\mu_p}=\frac{\sigma_{qi}}{\mu_q},\quad i\in I.
\end{equation}
where $(\sigma_{p1},...,\sigma_{pn})$ and $(\sigma_{q1},...,\sigma_{qn})$ are energies at $p$ and $q$, respectively.

\section{Proof of the a priori estimates and the degree counting theorems}

\noindent{\bf Proof of Theorem \ref{a-pri-main}:}

Let $u=(u_1,...,u_n)$ be a solution of (\ref{main-sys-2}). We set
\begin{equation}\label{v-u}
v_i=u_i-\log \int_M h_ie^{u_i}dV_g,\quad i=1,...,n,
\end{equation}
which immediately gives
\begin{equation}\label{int-e-1}
\int_Mh_ie^{v_i}dV_g=1,\quad i\in I.
\end{equation}
The equation for $v=(v_1,...,v_n)$ now becomes
\begin{equation}\label{liou-v}
\Delta_gv_i+\sum_{j\in I}\rho_j a_{ij}(h_je^{v_j}-1)=0,\quad \in I.
\end{equation}

To prove a priori estimate for $u$, we only need to establish
\begin{equation}\label{ap-v}
|v_i(x)|\le C,\quad i\in I,
\end{equation}
because with (\ref{ap-v}) we have
\begin{equation}\label{ap-u}
\log \int_M h_ie^{u_i}-C\le  u_i(x)\le \log \int_M h_ie^{u_i}+C.
\end{equation}
The fact that $u\in \, \hdot^{1,n}(M)$ implies that for each $i$, there exists $x_{0,i}\in M$ such that $u_i(x_{0,i})=0$. Hence by (\ref{ap-u}) we have
\begin{equation}\label{ap-u-2}
|\log \int_M h_ie^{u_i}| \le C,\quad i\in I.
\end{equation}
In view of (\ref{v-u}) and (\ref{ap-u-2}), the bound for $u$ is a direct consequence of the bound of $v$. Also  we only need to prove the upper bound for $v$, because the lower bound of $v$ can
be obtained from the upper bound of $v$ and standard Harnack inequality.
Therefore our goal is to prove
\begin{equation}\label{ub-v}
v_i(x)\le C,\quad i\in I.
\end{equation}

The proof of (\ref{ub-v}) is by contradiction. Suppose there exists a sequence $v^k$ to (\ref{liou-v}) that $\lim_{k\to \infty} \max_i\max_xv_i^k(x)\to \infty$. Then we consider two separate cases.

\medskip

\noindent {\emph{ Case one:} $\rho_i^k\to \rho_i>0$ as $k\to \infty$, for all $i\in I$}.

\medskip

The equation for $v^k$ is
\begin{equation}\label{e-vk}
\Delta_g v_i^k+\sum_{j\in I}\rho_j^k a_{ij}(h_je^{v_j^k}-1)=0,\quad i\in I.
\end{equation}
By an argument similar to a Brezis-Merle type lemma \cite{brezis-merle} it is easy to see that there are only finite blowup points: $\{p_1,...p_N\}$. Since $v_i^k$ is uniformly bounded above in any compact subset away from the blowup set, $v_i^k$ converges to
$\sum_{l=1}^Nm_{il}G(x,p_l)$ uniformly in compact sets away from $\{p_1,...,p_n\}$. Here we use the notation
$$\left\{\begin{array}{ll}
m_{il}=\sum_{j\in I}a_{ij}\sigma_{jl},\\
\\
\sigma_{il}=\lim_{k\to \infty} \frac{1}{2\pi}\int_{B(p_l,\delta)}\rho_j^kh_je^{v_j^k}dV_g,
\end{array}
\right.
$$
for some $\delta>0$, such that $B(p_l,2\delta)\cap B(p_s,\delta)=\emptyset$ for all $l\neq s$. To apply the local estimate we rewrite the equation for $v_i^k$ in local coordinates. For $p\in M$, let $y=(y^1,y^2)$ be the isothermal coordinates near $p$ such that $y_p(p)=(0,0)$ and $y_p$ depends smoothly on $p$. In this coordinates $ds^2$ has the form
$$e^{\phi(y_p)}[(dy^1)^2+(dy^2)^2], $$
where
$$\nabla \phi(0)=0,\quad \phi(0)=0.  $$
Also near $p$ we have
$$\Delta_{y_p}\phi=-2Ke^{\phi}, \quad \mbox{where } K \mbox{ is the Gauss curvature}. $$
When there is no ambiguity we write $y=y_p$ for simplicity. In local coordinates, the equation for $v_i^k$ can be written as
\begin{equation}\label{e-v-local}
-\Delta v_i^k=e^{\phi}\sum_{j=1}^na_{ij}\rho_j^k(h_je^{v_j^k}-1), \quad \mbox{in}\quad B(0,\delta),\quad i\in I.
\end{equation}

Let $f_i^k$ solve
$$-\Delta f_i^k=-e^{\phi}\sum_{j\in I}\rho_j^ka_{ij},\quad \mbox{in}\quad B(0,\delta),\quad i\in I, $$
and $f_i^k(0)=|\nabla f_i^k(0)|=0$. Set $\tilde v_i^k=v_i^k-f_i^k$ and
$$H_i^k=e^{\phi}\rho_i^ke^{f_i^k}h_i, $$
then the equation for $\tilde v_i^k$ becomes
\begin{equation}\label{t-vik}
-\Delta \tilde v_i^k=\sum_{j\in I}a_{ij}H_j^ke^{\tilde v_j^k},\quad \mbox{in}\quad B(0,\delta).
\end{equation}
Here we observe that
$$\int_{B(0,\delta)}H_i^ke^{\tilde v_i^k}dx=\int_{B(0,\delta)}\rho_i^kh_ie^{v_i^k}dV_g. $$
Since $v_i^k$ tends to $-\infty$ in $M\setminus \cup_{j=1}^NB(p_j,\delta)$, we have

\begin{equation}\label{fin-osi}
|\tilde v_i^k(x)-\tilde v_i^k(y)|\le C,\quad \forall x,y\in M\setminus \cup_{j=1}^NB(p_j,\delta/2),\quad i\in I.
\end{equation}

By Lemma \ref{mi-big-1} and the proof of Lemma \ref{local-est-1} it is easy to see than
\begin{equation}\label{lit-e}
\int_{M\setminus \cup_{j=1}^NB(p_j,\delta)}h_ie^{v_i^k}dV_g\to 0,\quad i\in I.
\end{equation}
and
\begin{equation}\label{quan-e}
\lim_{k\to \infty}\int_{B(p_l,\delta)}\rho_i^kh_ie^{v_i^k}dV_g/\mu_{p_l}=\lim_{k\to \infty}\int_{B(p_m,\delta)}\rho_i^kh_ie^{v_i^k}dV_g/\mu_{p_m}
\end{equation}
for $i\in I$ and any pair of $l,m$ between $1$ and $N$.

If we use $\mu_{p_l}$ to represent the possible strength of the singular source at each $p_l$, by  (\ref{same-e-2}) we have, for each $i\in I$,
$$\frac{\sigma_{i,1}}{\mu_1}=\frac{\sigma_{i,2}}{\mu_2}=...=\frac{\sigma_{i,N}}{\mu_N}, $$
and
$$2\pi(\sigma_{i,1}+\sigma_{i,2}+...+\sigma_{i,N})=\rho_i. $$
Thus
$$\sigma_{i,l}=\frac{\rho_i \mu_i}{2\pi\sum_{s=1}^n\mu_s},\quad i\in I,\quad l=1,...,N. $$
For each $l$, the Pohozaev identity for $(\sigma_{1,l},...,\sigma_{n,l})$ can be written as
$$\sum_{i,j\in I}a_{ij}\frac{\sigma_{i,l}}{\mu_i}\frac{\sigma_{j,l}}{\mu_j}=4\sum_{i\in I}\frac{\sigma_{i,l}}{\mu_i}. $$
Thus if blowup does happen, $(\rho_1,...,\rho_n)$ satisfies
\begin{equation}\label{pi-rho}
\sum_{i,j\in I}a_{ij}\rho_i\rho_j=8\pi\sum_{l=1}^N\mu_l\sum_{i\in I}\rho_i.
\end{equation}
Thus if $\rho$ is not  on critical hyper-surfaces $\Gamma_k$, the a priori estimate holds in this case.

\medskip

\noindent{\emph{Case two}:} Some of $\rho_i^k$ tend to $0$. Without loss of generality we assume that
 $\lim_{k\to \infty} \rho_i^k=\rho_i>0$, \quad $i\in I_1:=\{1,...,l\}$, $\lim_{k\to \infty} \rho_i^k=0$ for $i>l$.

 Let $M_k=\max\{v_1^k,....,v_l^k\}$ and $\bar M_k=\max\{v_{l+1}^k,...,v_n^k\}$. We first show that
 \begin{equation}\label{mk-bar}
 \bar M_k-M_k\le C.
 \end{equation}
If (\ref{mk-bar}) is not true, we have $\bar M_k-M_k\to \infty$, then we let
$$V_i^k(y)=v_i^k(e^{-\bar M_k/2}y+p_k)-\bar M_k $$
where $p_k$ is where $\bar M_k$ is attained: $v_{i_0}^k(p_k)=\bar M_k$. Clearly $i_0>l$. Thanks to the fact that $V_i^k\to -\infty$ for $i\le l$ and $\rho_i^k\to 0$ for $i>l$, $V_{i_0}^k$ converges uniformly to
$$\left\{\begin{array}{ll}
-\Delta V_{i_0}=0,\quad \mbox{in }\quad \mathbb R^2, \\
\\
V_{i_0}(0)=0.
\end{array}
\right.
$$
The fact that $V_{i_0}\equiv 0$ in $\mathbb R^2$ contradicts the finite energy of the component $i_0$. Thus (\ref{mk-bar}) is established.

We use the same notation as in Case one. Let $p_1$,...,$p_N$ be blowup points for $v_i^k$. The around each blowup point, say, $p_1$, the equation for $v^k$ can be written in local coordinates as (\ref{t-vik}) with $\tilde v_i^k$ and $H_i^k$ defined as in case one. Without loss of generality we assume that $\rho_i^k>0$ for all $k$ and $l+1\le i\le L$ and $\rho_i^k=0$ for all $k$ and $i>L$. Then we observe from the definition of $H_i^k$ and $H_i^k\to 0$ for $l+1\le i\le L$ and $H_i^k=0$ for $i>L$.

To reduce case two to case one, we need to adjust the terms involving vanishing $H_i^k$s. To do this we let $\hat f_i^k$ as
$$\left\{\begin{array}{ll}
-\Delta \hat f_i^k=\sum_{j=L+1}^na_{ij}e^{\tilde v_j^k-M_k},\quad \mbox{in}\quad B(0,\delta),\\
\\
\hat f_i^k(x)=0,\quad \mbox{on}\quad \partial B(0,\delta).
\end{array}
\right.
$$
Since $\max_i v_i^k-M_k$ is bounded above for all $i$, we have
$$\|\hat f_i^k\|_{C^1}\le C $$
for some $C$ independent of $k$. Now we define
$$\hat v_i^k=\left\{\begin{array}{ll}
\tilde v_i^k+\hat f_i^k,\quad i=1,...,l,\\
\tilde v_i^k+\log \rho_i^k+\hat f_i^k,\quad l+1\le i\le L,\\
\tilde v_i^k-M_k+\hat f_i^k,\quad L+1\le i\le n.
\end{array}
\right.
$$
and
$$\hat H_i^k=\left\{\begin{array}{ll}
H_i^ke^{-\hat f_i^k},\quad 1\le i\le l,\\
\frac{H_i^k}{\rho_i^k}e^{-\hat f_i^k}=e^{\phi+f_i^k-\hat f_i^k}h_i,\quad l+1\le i\le L,\\
e^{\hat f_i^k},\quad L+1\le i\le n.
\end{array}
\right.
$$
The definition of $\hat H_i^k$ immediately gives
$$\frac{1}c\le \hat H_i^k\le c,\quad \mbox{in}\quad B(0,\delta) $$
for some $c>0$ independent of $k$. Next
the equation for $\tilde v_i^k$ is
$$-\Delta \hat v_i^k=\sum_{j\in I}a_{ij}\hat H_j^k e^{\hat v_j^k},\quad \mbox{ in }\quad B(0,\delta),\quad i\in I.
$$

It is easy to see that $\max \hat v_i^k-M_k\to -\infty$ for $l+1\le i\le n$. Therefore case two is reduced to case one, which gives
$$\sigma_{il}/\mu_l=\sigma_{im}/\mu_m,\quad \forall l,m\in \{1,...,N\},\quad 1\le i\le l, $$
and $\sigma_{im}=0$ for all $i>l$ and all $m\in \{1,...,N\}$. Then as in case one if $(\rho_1,...,\rho_l,0,..,0)$ is not on any critical hyper-surfaces, the a priori estimate holds. Theorem \ref{a-pri-main} is established. $\Box$
\medskip

\noindent{\bf Proof of Theorem \ref{main-thm}:}
The main idea of the proof of the degree counting theorem is to reduce the whole system to the single equation.

\noindent{\emph{Case one:} At least one of $a_{ii}>0$}. We may assume $a_{11}>0$.  Thanks to Theorem \ref{a-pri-main} the Leray-Schauder degree of (\ref{main-sys-2}) for $\rho\in \mathcal{O}_{k}$ is equal to the degree for the following specific system corresponding to $(\rho_1,0,...,0)$:
\begin{equation}\label{main-sys-3}
\left\{\begin{array}{ll}
\Delta_gu_1+\rho_1 a_{11} (\frac{h_1e^{u_1}}{\int_M h_1e^{u_1}dV_g}-1)=0,\\
\Delta_gu_j+\rho_1a_{j1}(\frac{h_1e^{u_1}}{\int_M h_1e^{u_1}dV_g}-1)=0,\quad \mbox{for } j\ge 2,
\end{array}
\right.
\end{equation}
where $\rho_1$ satisfies
$$8\pi n_k<a_{11}\rho_1<8\pi n_{k+1}. $$
It is easy to see that $(\rho_1,0,...,0)\in \mathcal{O}_k$, using the degree counting formula of Chen-Lin \cite{chenlin3} for the single equation, we obtain the desired formula.

\noindent{\emph{Case two}: $a_{ii}=0$ for all $i\in I$}.

Using $a_{12}>0$, we reduce the degree counting formula for $\rho\in \mathcal{O}_k$ to the following system:
\begin{equation}\label{main-sys-4}
\left\{\begin{array}{ll}
\Delta_g u_1+a_{12}\rho_2(\frac{h_2e^{u_2}}{\int_M h_2e^{u_2}dV_g}-1)=0,\\
\Delta_g u_2+a_{12}\rho_1(\frac{h_1e^{u_1}}{\int_M h_1e^{u_1}dV_g}-1)=0,\\
\Delta_g u_i+\rho_1 a_{i1}(\frac{h_1e^{u_1}}{\int_M h_1e^{u_1}dV_g})+\rho_2 a_{12}(\frac{h_2e^{u_2}}{\int_M h_2e^{u_2}dV_g}-1)=0, \quad i\ge 3.
\end{array}
\right.
\end{equation}
where $\rho_1$, $\rho_2$ satisfy
$$8\pi n_k(\rho_1+\rho_2)<2 a_{12}\rho_1\rho_2<8\pi n_{k+1} (\rho_1+\rho_2). $$
It is easy to see that $(\rho_1,\rho_2,0,...,0)\in \mathcal{O}_k$. Now we consider the special case $\rho_1=\rho_2$, $h_1=h_2=h$. In this case a simple application of the maximum principle gives $u_1=u_2+C$, since they both have average equal to $0$, we have $u_1=u_2$. Then the first two equations in (\ref{main-sys-4}) turn out to be
$$\Delta_g u+a_{12}\rho(\frac{he^u}{\int_M he^u dV_g}-1)=0, $$
where $\rho\in (8\pi n_k,8\pi n_{k+1})$. Again the degree counting formula of Chen-Lin \cite{chenlin3} for the single equation gives the desired formula. Theorem \ref{main-thm} is established. $\Box$

\begin{rem}
The proof of Theorem \ref{thm-local} requires that there is no blowup point on $\partial \Omega$. Since all the singular sources are in the interior of $\Omega$, a standard moving plane argument can be employed to prove this fact. The interested readers may read into \cite{linzhang2} for the detail of the proof. Then the remaining part is similar to the proof of Theorem \ref{main-thm}.
\end{rem}

Finally we prove Theorem \ref{thm-spe}: Since the genus of the torus $M$ is $1$, $\chi(M)=0$ and the generating function is
\begin{eqnarray*}
g(x)=&\Pi_{p=1}^N\frac{1-x^{\mu_p}}{1-x}=\Pi_{p=1}^N(1+x+x^2+...+x^{\gamma_p})\\
=&1+b_1x+b_2x^2+...+b_kx^k+...+x^m.
\end{eqnarray*}
where $m=\sum_{p}\gamma_p$.
Let $$\rho_i=(\sum_{j\in I} a^{ij})4\pi \sum_{p=1}^N \gamma_p, $$
it is easy to see that
$$8\pi n_k\sum_i\rho_i<\sum_{ij}a_{ij}\rho_i\rho_j<8\pi n_{k+1}\sum_i \rho_i$$
for $n_k=(m-1)/2$ and $n_{k+1}=(m+1)/2$. Thus the Leray-Schauder degree $d_{\rho}$ can be computed as
$$d_{\rho}=\sum_{l=0}^{(m-1)/2}b_l. $$
Using $b_{m-l}=b_l$ for $l=0,1,..,m$ we further write $d_{\rho}$ as
$$d_{\rho}=\frac 12 \sum_{l=1}^m b_l=\frac{g(1)}2=\frac{\Pi_{p=1}^N(1+\gamma_p)}2. $$
Theorem \ref{thm-spe} is established. $\Box $

\end{document}